\pdfoutput=1
\RequirePackage{ifpdf}
\ifpdf 
\documentclass[pdftex]{sigma}
\else
\documentclass{sigma}
\fi

\usepackage{bbm}
\usepackage{shadow}
\usepackage[all]{xy}
\usepackage{tikz-cd}

\newcommand{\R}{\mathbb{R}}
\newcommand{\RP}{\mathbb{RP}}
\newcommand{\CP}{\mathbb{CP}}
\newcommand{\HP}{\mathbb{HP}}

\newcommand{\Z}{\mathbb{Z}}
\newcommand{\ZZ}{\mathbb{Z}}
\newcommand{\Q}{\mathbb{Q}}
\newcommand{\QQ}{\mathbb{Q}}

\newcommand{\C}{\mathbb{C}}

\newcommand{\HH}{\mathbb{H}}

\DeclareMathOperator{\Gr}{Gr}

\DeclareMathOperator{\Map}{Map}

\DeclareMathOperator{\End}{End}
\DeclareMathOperator{\Spin}{Spin}

\newcommand{\Cl}{\mathcal{C}{\rm l}}
\newcommand{\CCl}{\mathbb{C}{\rm l}}

\numberwithin{equation}{section}

\newtheorem{Theorem}{Theorem}[section]
\newtheorem*{Theorem*}{Theorem}
\newtheorem{Corollary}[Theorem]{Corollary}
\newtheorem{Lemma}[Theorem]{Lemma}
\newtheorem{Proposition}[Theorem]{Proposition}
 { \theoremstyle{definition}

\newtheorem{Example}[Theorem]{Example}
\newtheorem{Remark}[Theorem]{Remark} }

\begin{document}
\allowdisplaybreaks

\newcommand{\arXivNumber}{2207.12946}

\renewcommand{\PaperNumber}{093}

\FirstPageHeading

\ShortArticleName{Topology of Almost Complex Structures on Six-Manifolds}

\ArticleName{Topology of Almost Complex Structures\\ on Six-Manifolds}

\Author{Gustavo GRANJA~$^{\rm a}$ and Aleksandar MILIVOJEVI\'{C}~$^{\rm b}$}

\AuthorNameForHeading{G.~Granja and A.~Milivojevi\'{c}}

\Address{$^{\rm a)}$~Center for Mathematical Analysis, Geometry and Dynamical Systems, Instituto Superior\\
\hphantom{$^{\rm a)}$}~T\'ecnico, Universidade de Lisboa, Av. Rovisco Pais, 1049-001 Lisboa, Portugal}
\EmailD{\href{mailto:gustavo.granja@tecnico.ulisboa.pt}{gustavo.granja@tecnico.ulisboa.pt}}
\URLaddressD{\url{http://www.math.tecnico.ulisboa.pt/~ggranja/}}

\Address{$^{\rm b)}$~Max Planck Institute for Mathematics, Vivatsgasse 7, 53111 Bonn, Germany}
\EmailD{\href{mailto:milivojevic@mpim-bonn.mpg.de}{milivojevic@mpim-bonn.mpg.de}}
\URLaddressD{\url{https://guests.mpim-bonn.mpg.de/milivojevic/}}

\ArticleDates{Received August 08, 2022, in final form November 20, 2022; Published online December 02, 2022}

\Abstract{We study the space of (orthogonal) almost complex structures on closed six-dimensional manifolds as the space of sections of the twistor space for a given metric. For a connected six-manifold with vanishing first Betti number, we express the space of almost complex structures as a quotient of the space of sections of a seven-sphere bundle over the manifold by a circle action, and then use this description to compute the rational homotopy theoretic minimal model of the components that satisfy a certain Chern number condition. We further obtain a formula for the homological intersection number of two sections of the twistor space in terms of the Chern classes of the corresponding almost complex structures.}


\Classification{32Q60; 53C27; 53C28; 55P62}

\section{Introduction}

A major open problem in differential geometry is to determine whether a closed almost complex manifold of dimension at least six always admits an integrable complex structure. By a~celebrated theorem of Newlander and Nirenberg, an almost complex structure is integrable if and only if it satisfies a certain system of first order partial differential equations codified by the vanishing of the Nijenhuis tensor. An understanding of the topology of the space of all almost complex structures may be useful in the search for those which are integrable. We take this as motivation for our study of the topology of the space of almost complex structures on a~six-manifold. The present work is a continuation of~\cite{FGM21} which focused on the six-sphere and touched upon the case of almost complex structures on six-manifolds with vanishing first Chern class.\looseness=-1

Let $M$ denote an oriented six-manifold and $\mathcal J(M)$ be the space of
all almost complex structures on $M$ inducing the given orientation.
Fixing a Riemannian metric on $M$, the inclusion $J(M) \hookrightarrow \mathcal J(M)$ of the subspace of almost complex structures which are orthogonal with respect to the metric is a homotopy equivalence, allowing us to analyze $\mathcal J(M)$ via a space of sections of a fiber bundle with compact fibers. Generally, over an oriented Riemannian $2n$-manifold $M$ one can consider the ${\rm SO}(2n)/{\rm U}(n)$ bundle
$Z_+(M) \to M$ whose fiber over any $x \in M$ is the set of (linear) orthogonal complex structures on $T_x M$ compatible with the orientation on $M$; the space $Z_+(M)$ is known as the \emph{$($positive$)$ twistor space} of $M$. Sections of this bundle correspond to orthogonal almost complex structures on $M$ compatible with the orientation.

The twistor space construction can be done in any dimension, but six-manifolds enjoy the special property that the fiber of the twistor space is likewise six-dimensional. Hence two orthogonal almost complex structures on a closed Riemannian manifold $M$ give embeddings of the manifold into its twistor space which generically intersect at a finite number of points. One of our main results (Theorem~\ref{intersection}) computes this intersection number in terms of the Chern classes of the two almost complex structures. Namely, denoting by $c_i$ and $c_i'$ the respective Chern classes, the homological intersection number is given by
\[ \int_M \tfrac{1}{8}\big( c_1^3 + c_1^2 c_1' - c_1c_1'^2 - c_1'^3 \big) + \tfrac{1}{2} \big( c_1 c_2' + c_1' c_2' \big) - c_3.\]
In particular, the homological self-intersection of an almost complex structure is given by $\int_M c_1c_2 - c_3$. Employing a theorem of Michelsohn and Salamon, we observe that negative self-intersection restricts the possible deformations of an orthogonal complex structure (Corollary~\ref{deform}).

In order to study the homotopy type of $J(M)$, in Section~\ref{section3} we use the existence of ${\rm spin}^{\rm c}$ structures on oriented four-manifolds and almost complex six-manifolds to describe the integral cohomology of their twistor spaces. This will include both the positive and the \emph{negative} twistor space $Z_-(M)$, i.e., the total space of the bundle of linear complex structures inducing the opposite orientation on the tangent space. Presumably these results are well known, but we were not able to find them in the literature in this generality.

We take particular care with orientations throughout, due to varying historical conventions (e.g., the original~\cite{AHS78}, and~\cite{Hi81}, refer to the negative twistor space as the ``twistor space'', while this terminology is reserved for the positive twistor space for instance in~\cite{LM}\footnote{Despite this, the statement of~\cite[Theorem~4.1]{AHS78} in~\cite[Theorem~IV.9.14]{LM} uses the convention of~\cite{AHS78}.} and~\cite{Sa84}), and to the authors' difficulty in verifying some formulas in~\cite[p.~135]{Hi81}; see Remark~\ref{hitchin}.

The twistor spaces $Z_+(M)$ and $Z_-(M)$ (whose homotopy types generally differ when the dimension of $M$ is divisible by four) carry two natural almost complex structures, differing by the choice of induced orientation on the fiber (see, e.g.,~\cite[Proposition~3.1]{Sa84}), known as the \emph{Atiyah--Hitchin--Singer} and \emph{Eells--Salamon} almost complex structures. Even though the standard complex structure and its negative on the fibers ${\rm SO}(4)/{\rm U}(2) \cong \CP^1$ and ${{\rm SO}(6)/{\rm U}(3) \cong \CP^3}$ are biholomorphic via complex conjugation, the Eells--Salamon almost complex structure is never integrable~\cite[Proposition~3.4]{Sa84}, in contrast to the Atiyah--Hitchin--Singer structure (this phenomenon is not unique to dimensions four and six).

Our second main result is a description of the rational homotopy type of (a given component of) the space of almost complex structures on six-manifolds, under the additional assumption that $b_1(M) = 0$ and $\int_M c_1c_2 - c_3 \neq 0$ (Theorem~\ref{Qhomotopy6}). For this we use an expression of $J(M)$ as the quotient of the space of sections of an $S^7$-bundle by an $S^1$-action (Theorem~\ref{fibration}), which one can compare with the result obtained for $c_1=0$ in~\cite[Proposition~5.1]{FGM21}.
Throughout we provide examples for all of the mentioned results.

\subsection*{Notation and conventions}

 The projectivization of a complex vector bundle $E$ is denoted by $\mathbb{P}(E)$. We denote the Chern classes of a complex vector bundle by~$c_i$, and the Pontryagin classes of a real vector bundle by~$p_i$; if the vector bundle is complex, by its Pontryagin classes we mean those of the underlying real bundle. The Euler class is denoted by~$\mathsf{e}$. The Euler characteristic of a space is denoted by~$\chi$, and~$\sigma$ is the signature of an oriented closed manifold. All manifolds are assumed connected unless otherwise stated.

We denote by $\Gamma(E)$ the space of sections of a bundle $E \to M$ (where $E$ is not necessarily a~vector bundle). If $s \colon M \to E$ is a section, $\Gamma(E)_s$ denotes the connected component of~$s$ in~$\Gamma(E)$.
For $E \to M$ a vector bundle with a metric, we write $S(E) \to M$ for the associated sphere bundle. We denote by $J(M)$ the connected component of a specific point in $\Gamma(Z_+(M))$, which will be clear from context. For a vector space~$V$ with an inner product, $J(V)$ will denote the space of all orthogonal (linear) complex structures on $V$. If $V$ is oriented, $J_+(V)$ (respectively~$J_-(V)$) denotes the connected component of~$J(V)$ consisting of elements which induce the given orientation on $V$ (respectively the opposite orientation). By $\Map(X,Y)$ we denote the space of unbased maps $X \to Y$ with the compact-open topology.

Cohomology is singular cohomology with integral coefficients unless another choice of coefficients is indicated. The evaluation of a cohomology class $\alpha$ on a homology class $A$ is denoted by $\langle \alpha, A \rangle$ or by $\int_A \alpha$. In the context of rational homotopy theoretic models, $\Lambda(x_i)$ denotes the free graded-commutative algebra on generators $x_i$.

We will make frequent use of the projective bundle formula, which states that for a complex vector bundle $E \to M$ of complex rank $k$, the cohomology of the projectivized bundle $\mathbb{P}(E) \xrightarrow{p} M$ is the free $H^*(M)$-algebra on one generator $x$ subject only to the relation $x^k + c_1(E)x^{k-1} + \cdots + c_{k-1}(E)x + c_k(E) = 0$, where $x$ is the first Chern class of the line bundle $\mathcal{O}_{E^*}(1)$ (see~\cite[Definition~15.13]{D12}). This line bundle is dual
to the tautological (Hopf) line subbundle of $p^*(E)$ whose fiber over $\ell \in \mathbb{P}(E_x)$ is $\ell \subset E_x$. It restricts to $\mathcal{O}(1)$ on each fiber $\CP^{k-1}$.

For a real vector space $V$ with an endomorphism $J$ squaring to $-\mathrm{Id}$, we denote by $\Lambda_\C^*(V)$ the (complex) exterior algebra of the complex vector space $(V,J)$. For an almost complex manifold $(M,J)$, $\Lambda_\C^\ast (TM)$ is the (complex) exterior bundle associated to $(TM,J)$.

\section{Preliminaries}\label{section2}
In this section we review some of the models for the space of orthogonal complex structures on a Euclidean space following~\cite[Section~IV.9]{LM}, to which we refer for further details.

Let $(V, (\cdot, \cdot ))$ be a $2n$-dimensional real inner product space. We will write \[J(V)=\big\{ J \in \End(V) \colon J^2=-I, \, J \text{ is orthogonal} \big\}\]
for the space of orthogonal complex structures on $V$. This space naturally has the structure of a complex algebraic variety, given by the isomorphism
\[ J(V) \xrightarrow{\phi} \Gr_n^{\rm Iso}(V\otimes \C) \]
with the Grassmannian of $n$-dimensional complex subspaces of $V\otimes \C$, which are isotropic with respect to the complex bilinear form on $V\otimes \C$ obtained from the inner product on $V$ by linear extension. The map $\phi$ is defined by
\[ \phi(J) = \{ v\otimes 1 + Jv \otimes {\rm i} \colon v \in V\} \]
(it assigns to $J$ the plane $V_J^{0,1} \subset V\otimes \C$, which is the $-{\rm i}$ eigenspace of the complex linear extension of $J$ to $V\otimes \C$).

There is a tautologous complex vector bundle on $J(V)$ defined by
\[J(V) \times V \xrightarrow{\pi_1} J(V),\]
where the fiber over $J\in J(V)$ is given the complex structure determined by $J$.
The isomorphism $\phi$ is covered by a bundle map to the dual (or conjugate) tautological bundle over the Grassmannian.

\begin{Lemma}\label{positchern}
Consider the universal bundle
\[
J\big(\R^{2n}\big)={\rm O}(2n)/{\rm U}(n) \xrightarrow{\iota} {\mathcal B}{\rm U}(n) \to {\mathcal B}{\rm O}(2n).
\]
Then
\begin{enumerate}\itemsep=0pt
\item[$(i)$] $\iota$ classifies the tautologous complex vector bundle on $J\big(\R^{2n}\big)$.
\item[$(ii)$] $\iota^*(c_1)$ generates a subgroup of index $2$ in the infinite cyclic second cohomology group of
each of the two components of $J\big(\R^{2n}\big)$.
\item[$(iii)$] $\iota^*(c_1)$ is positive. In particular, $\iota^*(c_1)^{\frac{n^2-n}{2}}$ evaluates on the orientation class of each of the two components of $J\big(\R^{2n}\big)$ to a positive integer.
\end{enumerate}
\end{Lemma}
\begin{proof}
$(i)$ The universal bundle can be realized by taking ${\mathcal B}{\rm U}(n)={\mathcal E}{\rm O}(2n)/{\rm U}(n)$. The universal bundle over ${\mathcal B}{\rm U}(n)$ is then given by ${\mathcal E}{\rm O}(2n)\times_{{\rm U}(n)} \R^{2n} \to {\mathcal E}{\rm O}(2n)/{\rm U}(n)$, and pulls back along the inclusion of the fiber $\iota$ to ${\rm O}(2n) \times_{{\rm U}(n)} \R^{2n} \to {\rm O}(2n)/{\rm U}(n)$, which is the homogeneous expression of the tautologous bundle.

$(ii)$ The inclusion-induced map ${\mathcal B}{\rm U}(n) \to {\mathcal B}{\rm O}(2n)$ factors through the
double cover ${\mathcal B}{\rm SO}(2n) \allowbreak \to {\mathcal B}{\rm O}(2n)$, and the two components of ${\rm O}(2n)/{\rm U}(n)$ are,
respectively, the fiber ${\rm SO}(2n)/{\rm U}(n)$ of ${\mathcal B}{\rm U}(n) \to {\mathcal B}{\rm SO}(2n)$ and the image
of this fiber under a lift to ${\mathcal B}{\rm U}(n)$ of the nontrivial deck transformation of
${\mathcal B}{\rm SO}(2n)$. Therefore it suffices to prove the statement for the component
${\rm SO}(2n)/{\rm U}(n)$. This follows from the Serre spectral sequence of the bundle
${\mathcal B}{\rm U}(n) \to {\mathcal B}{\rm SO}(2n)$ (see for instance~\cite[Theorem~III.6.11]{MiTo}).

$(iii)$ $\iota^*(c_1)$ is the pullback of $c_1$ of the dual tautological bundle under the embedding
\[J\big(\R^{2n}\big) \xrightarrow{\cong} \Gr_n^{\rm Iso}\big(\C^{2n}\big) \subset \Gr_n\big(\C^{2n}\big).\]
The $n^{\mathrm{th}}$ exterior power of the dual tautological bundle on $\Gr\big(\C^{2n}\big)$ gives the
Pl\"ucker embedding of the Grassmannian in projective space. Hence $\iota^*(c_1)$ is the first
Chern class of a line bundle on $J\big(\R^{2n}\big)$ whose sections embed $J\big(\R^{2n}\big)$ in projective space and is therefore positive.\qedhere
\end{proof}

The above identification of orthogonal complex structures with isotropic planes in
$V\otimes \C$ leads to another description of the space $J(V)$ which will play an
important role in this paper.

Let $\Cl(V)$ denote the Clifford algebra determined by the quadratic form $q(v)=-\|v\|^2$, and let $\CCl(V)=\Cl(V)\otimes_\R \C$
denote its complexification. Let $S_{\C}$ denote an irreducible ($2^n$-dimensional) $\CCl(V)$-module.

For each element, i.e., \emph{spinor}, $\sigma \in S_{\C}$ we have the map $V\otimes_\R \C \to S_{\C}$ given by right Clifford multiplication by $\sigma$. A spinor is said to be \emph{pure} if the associated kernel is half-dimensional. Let ${\mathcal P}S _{\C} \subset S_\C$ denote the subset of pure spinors. As the subspace $\ker (\cdot \sigma) \subset V\otimes \C$ is isotropic
there is a natural map from the projectivization of the set of pure spinors to the isotropic Grassmannian
\[
\mathbb{P}( {\mathcal P}S _\C ) \xrightarrow{\psi} \Gr^{\rm Iso}_{2n}(V\otimes \C)
\]
defined by $[\sigma] \mapsto \ker (\cdot \sigma)$. This is a map of algebraic varieties and is in fact an ${\rm SO}(2n)$-equivariant isomorphism~\cite[Proposition~IV.9.7]{LM}.

The space $J(V)$ has two connected components corresponding to the two possible orientations induced by the complex structure. A choice of orientation for $V$ leads to a decomposition
\[J(V) = J_+(V) \coprod J_-(V)\]
with $J_+(V)$ the component corresponding to the given orientation. In turn this decomposes the Grassmannian of isotropic subspaces into the components of positive and negative isotropic subspaces.

On the level of spinors this decomposition takes the following form.
The complex volume element in $\CCl(V)$ is the element
$\omega_\C = {\rm i}^n e_1 \cdots e_{2n}$, where
$\{e_i\}$ is any oriented orthonormal basis. As $\omega_{\C}^2=1$, the volume element decomposes
the module $S_C$ into a direct sum of $+1$ and $-1$ eigenspaces
\[
S_\C = S_\C^+ \oplus S_\C^-.
\]
The elements of the summands are called the positive and negative spinors, respectively.
Note that, since $\omega_\C$ anti-commutes with the action of $V \subset \CCl(V)$, both
$S_\C^+$ and $S_\C^-$ are ${\rm spin}^{\rm c}$ representations, called the positive and negative spinor representations, respectively.

Every pure spinor is necessarily either positive or negative and the isomorphisms $\phi$, $\psi$ give rise to isomorphisms of algebraic varieties
\[
\mathbb{P}\big({\mathcal P}S _\C^\pm\big) \cong J_\pm(V).
\]
In dimensions $2n=4,6$ all nonzero spinors are pure~\cite[Remark~IV.9.12]{LM} and the isomorphisms identify the spaces
of orthogonal complex structure with the disjoint union of the projectivization of positive
and negative spinors.

A point $J \in J(V)$ gives rise to a useful description (cf.~\cite[Chapter I, equation~(5.27)]{LM})
of the irreducible module $S_\C$.
We will write $\langle \cdot, \cdot \rangle$ for the Hermitian inner product on the complex vector space $(V,J)$ whose real part is $ (\cdot, \cdot )$.
Let $\Lambda_\C^*(V)$ denote the exterior algebra of the complex vector space $(V,J)$.
Contraction of a vector $v \in V$ with $\omega \in \Lambda_\C^k(V)$ is determined by the expression
\[
v \lrcorner (w_1 \wedge \cdots \wedge w_k) = \sum_{j=1}^k (-1)^{j+1}\langle v, w_j \rangle w_1 \wedge \cdots \wedge \widehat{w_j} \wedge \cdots \wedge w_k.
\]
The action $V \otimes_\R \Lambda_\C^*(V) \to \Lambda_\C^*(V)$ defined by
\begin{equation*}
v \cdot \omega = v \wedge \omega - v \lrcorner \, \omega
\end{equation*}
satisfies
\begin{equation*}
v\cdot(v \cdot \omega) = - \|v\|^2 \omega,
\end{equation*}
as one easily checks using an orthonormal basis for $V$ having a
real multiple of a nonzero $v$ as its first element. The complex linear extension of this action gives $\Lambda_\C^*(V)$ the structure of a $\CCl(V)$-module of complex dimension $2^n$. This is the dimension of the irreducible complex $\CCl(V)$-module, so $\Lambda_\C^*(V)$ is a convenient description of this module.

The following is an easy computation which we include for convenience.

\begin{Lemma}\label{plusmin}
Let $(V,J)$ be an oriented $2n$-dimensional Euclidean real vector space with an orthogonal complex structure $J$ compatible with the orientation and let $S_\C=\Lambda^*_\C(V)$ be the irreducible $\CCl(V)$-module described above. Then
\[
S_\C^+ = \bigoplus_{k \text{ even }} \Lambda_\C^k(V), \qquad S_\C^- = \bigoplus_{k \text{ odd }} \Lambda_\C^k(V).
\]
\end{Lemma}
\begin{proof}
Let $e_1, \ldots, e_n$ be an orthonormal basis for $(V,J)$. Then we can write the complex volume element as
\[
\omega_\C = {\rm i}^n e_1 Je_1 \cdots e_n Je_n .
\]
Consider the basis $e_{i_1} \wedge \cdots \wedge e_{i_k}$ for $\Lambda^k_{\C}$ with $1\leq i_1 < \cdots < i_k \leq n$. Let $1 \leq m \leq n$
and assume that $m \not \in \{i_1, \ldots, i_k \}$. Then
\[
\begin{aligned}
(e_m Je_m) \cdot e_{i_1} \wedge \cdots \wedge e_{i_k} & = e_m \cdot \big( Je_m \wedge e_{i_1} \wedge \cdots \wedge e_{i_k} \big) \\
& = -\big\langle e_m, Je_m \big\rangle e_{i_1} \wedge \cdots \wedge e_{i_k} = -i \ e_{i_1} \wedge \cdots \wedge e_{i_k}.
\end{aligned}
\]
On the other hand, if $m \in \{i_1, \ldots, i_k\}$, let $l$ be such that $m=i_l$. Then
\[
\begin{aligned}
(e_m Je_m) \cdot e_{i_1} \wedge \cdots \wedge e_{i_k} & = e_m \cdot \big( {-} (-1)^{l-1} \big\langle Je_m, e_{i_l} \big\rangle e_{i_1} \wedge \cdots \wedge \widehat{e_{i_l}} \wedge \cdots \wedge e_{i_k} \big)\\
& = e_m \wedge\big( (-1)^l (-{\rm i}) e_{i_1} \wedge \cdots \wedge \widehat{e_{i_l}} \wedge \cdots \wedge e_{i_k} \big) \\
& = (-1)^{l-1} (-1)^l(-{\rm i}) e_{i_1} \wedge \cdots \wedge e_{i_k} = {\rm i} e_{i_1} \wedge \cdots \wedge e_{i_k}.
\end{aligned}
\]
Hence
\[
\omega_\C \cdot e_{i_1} \wedge \cdots \wedge e_{i_k} = {\rm i}^n {\rm i}^k (-{\rm i})^{n-k} e_{i_1} \wedge \cdots \wedge e_{i_k} = (-1)^{2n-k} e_{i_1} \wedge \cdots \wedge e_{i_k},
\]
which completes the proof.
\end{proof}

\begin{Remark}
\label{lift}
For $(V,J)$ an oriented $2n$-dimensional Euclidean real vector space with an orthogonal complex structure $J$ compatible with the orientation, the line in $S^+_\C$ corresponding to $J\in J_+(V)$ is $\Lambda^0_\C(V)\subset S_\C^+$.
Indeed,
\[(v \otimes 1 + w \otimes {\rm i}) \cdot 1 = v+Jw = 0 \ \Leftrightarrow \ w=Jv. \]
\end{Remark}

\section{Twistor spaces of four- and six-manifolds}\label{section3}

In this section we describe the integral cohomology of the twistor space of (Riemannian) almost complex six-manifolds, along with the integral cohomology of the twistor space of arbitrary oriented (Riemannian) four-manifolds. Finally we compute the Chern classes of the Atiyah--Hitchin--Singer almost complex structure on the latter.

\subsection{Spinors on almost complex four- and six-manifolds}

Let $M$ be an oriented Riemannian manifold and let $Z(M)$ denote the space of orthogonal almost complex structures on $M$, called the \emph{twistor space}. The space
$Z(M)$ has two components corresponding to the structures which induce the given orientation on $M$ or the opposite orientation. Recall, we
denote these by $Z_{\pm}(M)$, and refer to them as the \emph{positive} and \emph{negative} twistor spaces.

The Atiyah--Hitchin--Singer almost complex structure on $Z(M)$ is defined as follows: the Levi-Civita connection induces a connection on the bundle $Z(M) \to M$ splitting $TZ(M)$ into a~direct sum of horizontal and vertical subspaces
\[
TZ(M) =T^{\rm h} Z(M) \oplus T^{\rm v}Z(M).
\]
The almost complex structure on $T^{\rm v} Z(M)$ is determined by the algebraic variety structure of the fibers $J(T_x M)$ of $Z(M) \to M$ explained in Section~\ref{section2}. At a point $J_x$ in the fiber
over $x\in M$ the complex structure on $T^{\rm h}_{J_x} Z(M) = T_x M$ is given\footnote{The Eells--Salamon almost complex structure is obtained replacing the vertical component of the Atiyah--Hitchin--Singer almost complex structure by its negative; see~\cite[Section~3]{Sa84} for details.} by $J_x$.

A ${\rm spin}^{\rm c}$ structure on $M$~\cite[Appendix~D]{LM} yields a complex spinor
bundle $S_\C(M)\to M$ together with an action
\[
\CCl(TM) \otimes S_\C(M) \to S_\C(M).
\]
The discussion in the previous section then provides a one-to-one correspondence between projectivized pure spinors on $M$ and orthogonal complex structures inducing a given orientation:
\[
Z_{\pm}(M) \cong \mathbb{P}\big( {\mathcal P}S _{\C}^{\pm}(M)\big).
\]
Now suppose $M$ is given an orthogonal almost complex structure $J$ compatible with the orientation. Then $J$ gives rise to a canonical ${\rm spin}^{\rm c}$ structure on $M$ via the group homomomorphism ${\rm U}(n) \xrightarrow{\kappa} \Spin^{\rm c}(2n)$, which is the unique
lift as a group homomorphism in the following diagram{\samepage
\[
\begin{tikzcd}
& \Spin^{\rm c}(2n) = \big(\Spin(2n) \times S^1\big)/\{\pm(1,1)\} \ar[d, "\pi \times \delta"] \\
{\rm U}(n) \ar[ur,"\kappa"] \ar[r,"\iota \times \det"] & {\rm SO}(2n) \times S^1.
\end{tikzcd}
\]
Here $\iota$ denotes the inclusion, $\pi\big(\big[h, {\rm e}^{{\rm i}\theta}\big]\big) = [h]$, and $\delta\big(\big[h,{\rm e}^{{\rm i}\theta}\big]\big)={\rm e}^{2{\rm i}\theta}$.}

An element $g\in {\rm U}(n)$ is diagonalized by some orthonormal basis $(e_1,Je_1, \ldots, e_n, Je_n)$. If
\[g=\mathrm{diag}\big({\rm e}^{{\rm i}\theta_1},\ldots, {\rm e}^{{\rm i}\theta_n}\big),\]
then
\[
\kappa(g)= \left[ \prod_{k=1}^n \left( \cos \frac{\theta_k} 2 + {\rm i} \sin \frac{\theta_k}{2} e_k J e_k\right), {\rm e}^{{\rm i} \frac{\sum \theta_k}{2}} \right]
\]
 (cf.~\cite[equation~(D.10)]{LM}). Using
the formulas for the action of $e_i Je_i$ in the proof of Lemma~\ref{plusmin}, we see that the spinor action of $\kappa(g)$
on $\Lambda^*_{\C}\big(\R^{2n}\big)$, where $\R^{2n}$ is equipped with the complex structure given by $\mathrm{diag}\left( \left(\begin{smallmatrix} 0 & -1 \\ 1 & 0 \end{smallmatrix}\right), \ldots, \left(\begin{smallmatrix} 0 & -1 \\ 1 & 0 \end{smallmatrix}\right) \right)$, agrees with the standard action of $g$ on this space. We obtain the following result:

\begin{Proposition}\label{identification}
Let $M$ be an oriented Riemannian manifold with an orthogonal almost complex structure $J$ compatible with the orientation. Then we have the following isomorphisms of smooth fiber bundles over $M$:
\begin{enumerate}\itemsep=0pt
\item[$1.$] If $\dim(M)=4$, then
\[
Z_+(M) \cong \mathbb{P}\big( \C \oplus \Lambda^2_\C(TM)\big), \qquad Z_{-}(M) \cong \mathbb{P}( TM).
\]
\item[$2.$] If $\dim(M)=6$, then
\[
Z_+(M) \cong \mathbb{P}\big( \C \oplus \Lambda^2_\C(TM)\big), \qquad Z_{-}(M) \cong \mathbb{P}\big( TM \oplus \Lambda_\C^3(TM)\big).
\]
\end{enumerate}
\end{Proposition}

\begin{Remark}
 The Atiyah--Hitchin--Singer and Eells--Salamon almost complex structures on (either) twistor space do not agree with the natural almost complex structure on the projectivization. Indeed, the projection map from the projectivization is complex linear whilst this is not the case for the twistor projection.
\end{Remark}

The projective bundle formula now gives the following formulas for the cohomology rings of the components of the twistor space in terms of the
Chern classes of $(TM, J)$.

\begin{Corollary}
\label{cohomtwis}
Let $M$ be an oriented Riemannian manifold with an orthogonal almost complex structure $J$ compatible with the orientation. We have the following isomorphisms as algebras over $H^*(M)$:
\begin{enumerate}\itemsep=0pt
\item[$1.$] If $\dim(M)=4$, then
\begin{align*}
&H^*(Z_+(M))\cong H^*(M)[x]/\big(x^2+c_1 x\big), \\
&H^*(Z_{-}(M))\cong H^*(M)[x]/\big(x^2 + c_1 x + c_2\big).
\end{align*}
\item[$2.$] If $\dim(M)=6$, then
\begin{align*}
&H^*(Z_+(M))\cong H^*(M)[x]/\big(x^4 +2c_1x^3 +\big(c_1^2+c_2\big) x^2 + (c_1c_2-c_3)x\big), \\
&H^*(Z_{-}(M))\cong H^*(M)[x]/\big(x^4 + 2c_1x^3+ \big(c_1^2+c_2\big)x^2 + (c_1c_2+c_3)x\big). \end{align*}
\end{enumerate}
\end{Corollary}

\begin{Remark}\quad
 \begin{itemize}\itemsep=0pt \item The formulas above in dimension four agree with~\cite[Theorem~11.2]{ESa85}, where it is established that over a K\"ahler surface $M$, the positive twistor space is the projectivization of $T^{2,0}M \oplus \C$ (note that $T^{2,0}M \cong \Lambda_{\C}^2(TM)$), and the negative twistor space is the projectivization of the holomorphic tangent bundle.

\item In~\cite[Section~6]{Ev14}, Evans describes the cohomology ring of the positive twistor space of a~(not necessarily almost complex) six-manifold, with complex coefficients. The description is in terms of the first Chern class of $T^{\rm h}(Z_+(M))$. Using Proposition~\ref{chernclasses}, one can check that Evans' description agrees with the one in Corollary~\ref{cohomtwis} in the case of almost complex manifolds. We remark that Evans' computation can quickly be reproduced using naturality in the pullback diagram
\[\begin{tikzcd}
Z_+(M) \arrow[d] \arrow[r] & {\mathcal B}{\rm U}(3) \arrow[d] \\
M \arrow[r] & {\mathcal B}{\rm SO}(6)
\end{tikzcd}\]
and expressing $H^*({\mathcal B}{\rm U}(3);\Q)$ as an $H^*({\mathcal B}{\rm SO}(6);\Q)$-algebra using the known behavior of the right-hand side vertical map on cohomology. \end{itemize}
\end{Remark}

\begin{Remark}\label{cp2bar} In order for the formulas of Corollary~\ref{cohomtwis}(1) to hold, it is crucial to have an almost complex structure \emph{inducing the given orientation} on the base manifold. For instance, the positive twistor space of $\overline{\CP^2}$ (i.e., the complex projective plane with the opposite orientation) with the Fubini--Study metric is the full flag variety ${\rm U}(3)/\left( {\rm U}(1)\times {\rm U}(1) \times {\rm U}(1) \right)$~\cite[p.~217, Example~1]{Sa84}, whose cohomology is given by
\[\ZZ[x,y,z]/(x+y+z, xy+xz+yz, xyz) \cong \ZZ[x,y]/\big(x^2+xy+y^2, x^2y+xy^2\big),\]
where $x$, $y$, $z$ are in degree two~\cite[Proposition~31.1]{B53}. This ring is not isomorphic to
\[\ZZ[x,y]/\big(y^3, x^2 + kyx\big)\]
for any integer $k$, so the cohomology of $Z_+\big(\overline{\CP^2}\big)$ can not be expressed as in the first item of Corollary~\ref{cohomtwis}(1). This corresponds to the fact that $\overline{\CP^2}$ does not admit an almost complex structure compatible with its orientation (which can of course be seen by other means as well, e.g., Hirzebruch's congruence $\chi + \sigma \equiv 0 \bmod 4$ for closed almost complex four-manifolds).
\end{Remark}

\subsection[Spin\^{}c(4)-bundles and twistor spaces of four-dimensional manifolds]{$\boldsymbol{\Spin^{\rm c}(4)}$-bundles and twistor spaces of four-dimensional manifolds}

Twistor spaces of oriented Riemannian four-manifolds give interesting examples of almost complex six-manifolds which we will consider in the following sections. In this subsection we will compute their integral cohomology as well as the Chern classes of
their Atiyah--Hitchin--Singer almost complex structures.

Our computation requires knowledge of the integral cohomology of ${\mathcal B}{\rm Spin}^{\rm c}(4)$. The cohomology of ${\mathcal B}{\rm Spin}^{\rm c}(n)$ is described in detail in~\cite{Du18}. We need only the simple case when $n=4$, so we include an elementary treatment of this case.

Our aim is to compute the cohomology ring of the ``universal twistor space" over ${\mathcal B}{\rm Spin}^{\rm c}(4)$, namely the bundle of orthogonal complex structures on the fibers of the $4$-plane bundle over ${\mathcal B}{\rm Spin}^{\rm c}(4)$ classified by the canonical map ${\mathcal B}{\rm Spin}^{\rm c}(4) \to {\mathcal B}{\rm SO}(4)$. Since every oriented four-manifold admits a ${\rm spin}^{\rm c}$ structure, naturality will
yield the cohomology ring of the twistor space.

We will use unit quaternions to describe the four-dimensional spin groups and orthogonal complex structures. We consider the quaternions $\HH$ as a vector space over $\C$ via
right multiplication by $\C \subset \HH$, i.e., we consider the identification of $\C^2$ with $\HH$ given by
\[ (z_1, z_2) \mapsto z_1 + {\rm j}z_2. \]
The standard complex structure on $\HH= \R^4$ is then right multiplication by ${\rm i}$,
and $(1,{\rm i},{\rm j},-{\rm k})$ is an oriented basis.

The unit quaternion $q= w_1 + {\rm j} w_2$ acts on $\HH$ by left multiplication as
\[ (w_1 + {\rm j}w_2) (z_1+{\rm j}z_2) = (w_1z_1 -\overline{w_2} z_2) + {\rm j}(w_2 z_1+z_2 \overline{w_1} ),\]
which corresponds via the identification with $\C^2$ to the matrix
\[
\left[ \begin{matrix}
w_1 & -\overline{w_2} \\
w_2 & \hphantom{-}\overline{w_1}
\end{matrix}
\right]
\]
with $|w_1|^2 +|w_2|^2=1$. This will be our identification of ${\rm Sp}(1)$ with ${\rm SU}(2)$.

The universal cover of ${\rm SO}(4)$ is modelled by the map ${\rm Sp}(1) \times {\rm Sp}(1) \to {\rm SO}(4)$
given by
\[ (q_1, q_2) \mapsto \left( v \mapsto q_1 v \overline{q_2} \right) \]
and we will identify
\[
{\rm SO}(4) = ({\rm Sp}(1) \times {\rm Sp}(1))/\{\pm (1,1)\}.
\]
In these terms, the subgroup ${\rm U}(2) \subset {\rm SO}(4)$ is
\[
{\rm U}(2) = \big({\rm Sp}(1) \times S^1\big)/\{\pm (1,1)\}.
\]

The group $\Spin^{\rm c}(4)$ is defined as
\[
\Spin^{\rm c}(4) = \big({\rm Sp}(1) \times {\rm Sp}(1) \times S^1\big)/\{\pm(1,1,1)\}.
\]
There are two canonical homomorphisms
\begin{equation}\label{canonichom}
{\rm SO}(4) \xleftarrow{\pi} \Spin^{\rm c}(4) \xrightarrow{\delta} S^1
\end{equation}
given by
\[
\pi\big( \big[q_1,q_2,{\rm e}^{{\rm i}\theta}\big]\big) = [q_1,q_2],
\qquad
\delta\big( \big[q_1,q_2, {\rm e}^{{\rm i}\theta}\big]\big) = {\rm e}^{2{\rm i}\theta}.
\]
The canonical map ${\rm U}(2) \xrightarrow{\kappa} \Spin^{\rm c}(4)$ is the unique lift as a group
homomorphism in the diagram
\[
\begin{tikzcd}
& \Spin^{\rm c}(4) \ar["\pi \times \delta",d] \\
{\rm U}(2) \ar["\kappa",ur] \ar["\iota \times \det", r] & {\rm SO}(4)\times S^1.
\end{tikzcd}
\]
In terms of the coordinates above, this is given by $\kappa\big(\big[q_1,{\rm e}^{{\rm i}\theta}\big]\big)=\big[q_1,{\rm e}^{{\rm i}\theta},{\rm e}^{-{\rm i}\theta}\big]$.

Let $x \in H^2\big({\mathcal B}S ^1\big)=H^2(\CP^\infty)$ denote the standard generator (the Chern class
of the dual tautological line bundle on $\CP^\infty$) and
\[
\alpha= ({\mathcal B}\delta)^*(x) \in H^2\big({\mathcal B}{\rm Spin}^{\rm c}(4)\big).
\]
We will write $p_1$ and $\mathsf{e}$ for the universal Pontryagin class and Euler class in
$H^*({\mathcal B}{\rm SO}(4))$, respectively.

\begin{Lemma}
\label{cohspin}
Let $\pi\colon \Spin^{\rm c}(4) \to {\rm SO}(4)$ be
as in \eqref{canonichom}. There are unique classes $S_1,S_2 \in H^4\big({\mathcal B}{\rm Spin}^{\rm c}(4)\big)$ such
that
\begin{align*}
&4S_1={\mathcal B}\pi^* p_1 - 2{\mathcal B}\pi^* \mathsf{e} - \alpha^2, \\
&4S_2={\mathcal B}\pi^* p_1 + 2{\mathcal B}\pi^* \mathsf{e} - \alpha^2.
\end{align*}
Moreover the cohomology ring of ${\mathcal B}{\rm Spin}^{\rm c}(4)$ is a polynomial algebra,
\[
H^*\big({\mathcal B}{\rm Spin}^{\rm c}(4)\big) \cong \Z[\alpha,S_1,S_2].
\]
\end{Lemma}
\begin{proof}
The short exact sequence of groups
\[
{\rm Sp}(1) \times {\rm Sp}(1) \xrightarrow{\iota} \Spin^{\rm c}(4) \xrightarrow{\delta} S^1,
\]
where $\iota(q_1,q_2)=[q_1,q_2,1]$, leads to a fiber sequence
\begin{equation}
\label{fs}
{\mathcal B}{\rm Sp}(1)\times {\mathcal B}{\rm Sp}(1) \xrightarrow{{\mathcal B}\iota} {\mathcal B}{\rm Spin}^{\rm c}(4) \xrightarrow{{\mathcal B}\delta} {\mathcal B}S ^1.
\end{equation}
Let $A,B \in H^4({\mathcal B}{\rm Sp}(1) \times {\mathcal B}{\rm Sp}(1))$ be the generators which map to
$x^2 \in H^4\big({\mathcal B}S ^1\big)$ under the maps induced by the
natural inclusions ${\rm e}^{{\rm i}\theta} \mapsto \big({\rm e}^{{\rm i}\theta},1\big)$ and ${\rm e}^{{\rm i}\theta} \mapsto
\big(1,{\rm e}^{{\rm i}\theta}\big)$ of $S^1$ in ${\rm Sp}(1) \times {\rm Sp}(1)$.
By the Serre spectral sequence of the fibration \eqref{fs}, any elements in
$H^4\big({\mathcal B}{\rm Spin}^{\rm c}(4)\big)$ mapping to $A$, $B$ under ${\mathcal B}\iota^*$ together with $\alpha$ will
freely generate the cohomology ring of ${\mathcal B}{\rm Spin}^{\rm c}(4)$.

The effect of the composition
\[ {\mathcal B}{\rm Sp}(1) \times {\mathcal B}{\rm Sp}(1) \xrightarrow{{\mathcal B}\iota} {\mathcal B}{\rm Spin}^{\rm c}(4) \xrightarrow{{\mathcal B}\pi} {\mathcal B}{\rm SO}(4) \]
on degree $4$ cohomology is determined by the homomorphism
\[
S^1 \times S^1 \xrightarrow{ \left[\begin{smallmatrix} \hphantom{-}1 & -1 \\ -1 & -1 \end{smallmatrix}\right] } S^1 \times S^1
\]
induced by $\pi \circ \iota$ on the standard maximal tori.
Writing $a,b \in H^2\big({\mathcal B}S ^1 \times {\mathcal B}S ^1\big)$ for the standard generators corresponding to the maximal torus of ${\rm SO}(4)$, and $x$, $y$ for
the corresponding generators for ${\rm Sp}(1)\times {\rm Sp}(1)$, we have
\[ a^2+b^2 \mapsto (x-y)^2+ (-x-y)^2= 2x^2 +2y^2, \qquad ab \mapsto -x^2+y^2, \]
and therefore $p_1,\mathsf{e} \in H^4({\mathcal B}{\rm SO}(4))$ map to $2A + 2B$ and $-A+B$ in $H^4({\mathcal B}{\rm Sp}(1)\times {\mathcal B}{\rm Sp}(1))$, respectively. Therefore
\[ {\mathcal B}\iota^* {\mathcal B}\pi^* (p_1 -2\mathsf{e}) = 4A, \qquad {\mathcal B}\iota^* {\mathcal B}\pi^* (p_1+2\mathsf{e}) = 4B.\]
It follows from the short exact sequence
\[
0 \to \Z\alpha^2 \to H^4\big({\mathcal B}{\rm Spin}^{\rm c}(4)\big) \xrightarrow{{\mathcal B}\iota^*} H^4({\mathcal B}{\rm Sp}(1)\times {\mathcal B}{\rm Sp}(1)) \to 0
\]
that there exists $\lambda \in \Z$ such that ${\mathcal B}\pi^*(p_1-2\mathsf{e}) + \lambda \alpha^2$
is (uniquely) divisible by $4$.

The composition
\[
\CP^2 \xrightarrow{\tau_{\CP^2}} {\mathcal B}{\rm U}(2) \xrightarrow{{\mathcal B}\kappa} {\mathcal B}{\rm Spin}^{\rm c}(4)
\xrightarrow{{\mathcal B}\delta} {\mathcal B}S ^1
\]
classifying the second exterior power of the tangent bundle of $\CP^2$ induces multiplication by $3$ on $H^2$, and hence $({\mathcal B}\kappa \circ \tau_{\CP^2})^*$ provides a $2$-local splitting of the map ${\mathcal B}\delta^*$ on $H^4$.

The pullback to $H^4\big(\CP^2\big)$ of
${\mathcal B}\pi^*(p_1 - 2\mathsf{e}) + \lambda \alpha^2 \in H^4\big({\mathcal B}{\rm Spin}^{\rm c}(4)\big)$ is
$(3-6+9\lambda)$ times the orientation class, and the smallest value of
$\lambda$ for which this number is a multiple of $4$ is $\lambda=-1$.

Since ${\mathcal B}\pi^*(p_1) - 2{\mathcal B}\pi^*(\mathsf{e}) - \alpha^2 \in H^4\big({\mathcal B}{\rm Spin}^{\rm c}(4)\big)$ is divisible
by $4$, we conclude that
\[
{\mathcal B}\iota^*\left( \frac{ {\mathcal B}\pi^*(p_1) - 2{\mathcal B}\pi^*(\mathsf{e}) - \alpha^2 }{4} \right)=A,
\]
and similarly
\[
{\mathcal B}\iota^*\left( \frac{ {\mathcal B}\pi^*(p_1) + 2{\mathcal B}\pi^*(\mathsf{e}) - \alpha^2 }{4} \right)=B,
\]
which completes the proof.
\end{proof}

For four-dimensional $M$, a choice of ${\rm spin}^{\rm c}$ structure gives us the following diagram of pullback squares:
\[
\begin{tikzcd}
Z(M) \ar[r] \ar[d] & Z\big({\mathcal B}{\rm Spin}^{\rm c}(4)\big) \ar[r] \ar[d] & Z({\mathcal B}{\rm SO}(4)) \ar[d] \\ M \ar[r] &
{\mathcal B}{\rm Spin}^{\rm c}(4) \ar[r,"{\mathcal B}\pi"] & {\mathcal B}{\rm SO}(4),
\end{tikzcd}
\]
where $Z({\mathcal B}G)$ denotes the bundle of orthogonal complex structures on the universal oriented $4$-plane bundle over ${\mathcal B}G$ for $G={\rm SO}(4)$ and its pullback for $G=\Spin^{\rm c}(4)$.

We will see that each of the two components of $Z\big({\mathcal B}{\rm Spin}^{\rm c}(4)\big)$ is the projectivization of a~complex plane bundle over ${\mathcal B}{\rm Spin}^{\rm c}(4)$. The projective bundle formula together with naturality will give us the following description of the cohomology ring of $Z(M)$ as an $H^*(M)$-algebra.

\begin{Proposition}
\label{cohtwist4}
Let $M$ be an oriented Riemannian four-manifold and $\alpha \in H^2(M)$ be an integral
lift of $w_2(M)$ $($classifying the complex line bundle associated to a ${\rm spin}^{\rm c}$-structure on $M)$. Then there is an isomorphism of $H^*(M)$-algebras
\begin{equation}
\label{cohtwis4}
H^*(Z_\pm(M)) \cong H^*(M)[x]/\left( x^2 + \alpha x - \frac{ p_1 \pm 2\mathsf{e} - \alpha^2}{4} \right).
\end{equation}
\end{Proposition}
\begin{proof}
Under our identification of $\R^4$ with $\HH$,
the orthogonal complex structures on $\R^4$ compatible with the orientation (respectively, opposite orientation) are given by right (respectively, left) multiplication by a unit imaginary quaternion. The element $[q_1,q_2] \in {\rm SO}(4)$ acts on $J_+\big(\R^4\big)$ by
\[
[q_1,q_2] u = q_2 u \overline q_2,
\]
and on $J_-\big(\R^4\big)$ by
\[
[q_1,q_2] u = q_1 u \overline q_1.
\]
That is, ${\rm SO}(4)$ acts on orthogonal complex structures via the two projections to ${\rm SO}(3)$.
The commutative diagrams
\[
\begin{tikzcd}
\Spin^{\rm c}(4) \arrow[d,"\pi"] \arrow[r,"\chi_+"] & {\rm U}(2) \arrow[d] & & \big[q_1,q_2,{\rm e}^{{\rm i}\theta}\big] \arrow[d, mapsto] \arrow[r, mapsto] & \big[q_2, {\rm e}^{-{\rm i}\theta}\big] \arrow[d, mapsto] \\
{\rm SO}(4) \arrow[r, "\pi_{\mathrm{right}}"] & {\rm SO}(3), & & \big[q_1,q_2\big] \arrow[r, mapsto] & \big[q_2\big].
\end{tikzcd}
\]
and
\[
\begin{tikzcd}
\Spin^{\rm c}(4) \arrow[d,"\pi"] \arrow[r,"\chi_-"] & {\rm U}(2) \arrow[d] & & \big[q_1,q_2,{\rm e}^{{\rm i}\theta}\big] \arrow[d, mapsto] \arrow[r, mapsto] & \big[q_1, {\rm e}^{-{\rm i}\theta}\big] \arrow[d, mapsto] \\
{\rm SO}(4) \arrow[r, "\pi_{\mathrm{left}}"] & {\rm SO}(3), & & \big[q_1,q_2\big] \arrow[r, mapsto] & \big[q_1\big].
\end{tikzcd}
\]
allow us to express each component of the bundle of orthogonal complex structures
associated to an oriented $4$-plane bundle with a ${\rm spin}^{\rm c}$ structure as the projectivization of a complex plane bundle. The formula \eqref{cohtwis4} will follow from the projective bundle formula once we compute the Chern classes of these plane bundles. By naturality it suffices to compute the images of~$c_1$,~$c_2$ under ${\mathcal B}\chi_{\pm}^*$.

As the composition $\Spin^{\rm c}(4) \xrightarrow{\chi_\pm} {\rm U}(2) \xrightarrow{\det} S^1 $
equals $\delta$, we have ${\mathcal B}\chi_{\pm}^*(c_1) = \alpha$.

We will now show that
\[ {\mathcal B}\chi_+^*(c_2) =-S_2, \qquad {\mathcal B}\chi_-^*(c_2)=-S_1.\]
In view of Lemma~\ref{cohspin}, this will complete the proof. The maps
\[
{\rm Sp}(1) \times {\rm Sp}(1) \xrightarrow{\iota} \Spin^{\rm c}(4) \xrightarrow{\chi_\pm} {\rm U}(2)
\]
are given respectively by $(q_1,q_2) \mapsto [q_2,1]$ and $(q_1,q_2) \mapsto [q_1,1]$ (i.e., they correspond respectively to
the inclusions of the right and left copies of ${\rm SU}(2)$ in ${\rm U}(2)$). Therefore, with respect to the standard basis for $H^4 ({\mathcal B}{\rm Sp}(1)\times {\mathcal B}{\rm Sp}(1))$ used in the proof of Lemma~\ref{cohspin}, we
have ${\mathcal B}(\iota \circ \chi_+)^*c_2 = -B$ and ${\mathcal B}(\iota\circ \chi_-)^*c_2 =-A$.

Hence
\[
{\mathcal B}\chi_+^*(c_2) = -S_2 + \lambda \alpha^2, \qquad {\mathcal B}\chi_-^*(c_2) = -S_2 + \mu \alpha^2, \qquad \text{for some } \lambda, \mu \in \Z.
\]
We can determine the coefficients $\lambda$ and $\mu$ by mapping to $H^4\big(\CP^2\big)$ under $\CP^2 \! \xrightarrow{{\mathcal B}\kappa \circ \tau_{\CP^2}} \! {\mathcal B}{\rm Spin}^{\rm c}(4)$.
The compositions
\[ {\rm U}(2) \xrightarrow{\kappa} \Spin^{\rm c}(4) \xrightarrow{\chi_{\pm}} {\rm U}(2) \]
are given by
\[
\big[q, {\rm e}^{{\rm i}\theta}\big] \mapsto \big[{\rm e}^{{\rm i}\theta},{\rm e}^{{\rm i}\theta}\big] \qquad \text{and} \qquad \big[q, {\rm e}^{{\rm i}\theta}\big] \mapsto \big[q,{\rm e}^{{\rm i}\theta}\big],
\]
respectively. The second map is the identity, while the first map is the representation ${\rm U}(2) \xrightarrow{ 1\oplus \det } {\rm U}(2)$ (cf.\ Proposition~\ref{identification}(1)).
It follows that the composition
\[ \CP^2 \to {\mathcal B}{\rm Spin}^{\rm c}(4) \xrightarrow{{\mathcal B}\chi_+} {\mathcal B}{\rm U}(2) \]
classifies the bundle $\C \oplus \Lambda^2 T\CP^2$, and
\[ \CP^2 \to {\mathcal B}{\rm Spin}^{\rm c}(4) \xrightarrow{{\mathcal B}\chi_-} {\mathcal B}{\rm U}(2) \]
classifies the tangent bundle of $\CP^2$. Thus $c_2 \in H^4({\mathcal B}{\rm U}(2))$ must go to $0 \in H^4\big(\CP^2\big)$ under the composition $\CP^2 \to {\mathcal B}{\rm Spin}^{\rm c}(4) \xrightarrow{{\mathcal B}\chi_+} {\mathcal B}{\rm U}(2)$, and to $3\in H^4\big(\CP^2\big)$
under the composition $\CP^2 \to {\mathcal B}{\rm Spin}^{\rm c}(4) \xrightarrow{{\mathcal B}\chi_-} {\mathcal B}{\rm U}(2)$, i.e.,
\begin{align*}
- (3 + 6 - 9)/4 + 9 \lambda = 0 \ \Leftrightarrow \ \lambda=0, \qquad\! \mathrm{and} \qquad\! {-}(3-6-9)/4 + 9\mu &= 3 \ \Leftrightarrow \ \mu = 0.\!\!\!\tag*{\qed}
 \end{align*}\renewcommand{\qed}{}
\end{proof}

Recall that $T^{\rm v} Z_{\pm}(M)$ and $T^{\rm h}Z_{\pm}(M)$ denote the vertical and horizontal subbundles of \linebreak $TZ_{\pm}(M)$ which are complex vector bundles via the Atiyah--Hitchin--Singer almost complex structure.

\begin{Proposition}\label{chernclasses}
Let $M$ be an oriented Riemannian four-manifold and $\alpha \in H^2(M)$ a choice of lift of $w_2(M)$. Then, in terms of the expression for $H^*(Z_\pm(M))$ in Proposition~{\rm \ref{cohtwist4}}, the total Chern class of the standard almost complex structure on $Z_\pm(M)$ is
\begin{equation}
\label{chernform}
1 + (4x+2\alpha) + (p_1 \pm 3\mathsf{e}) \pm (\alpha + 2x)\mathsf{e}.
\end{equation}
Moreover, $c_1(T^{\rm v}Z_\pm(M))=c_1\big(T^{\rm h}Z_\pm(M)\big)=\alpha + 2x$.
\end{Proposition}
\begin{proof}
It suffices to compute the total Chern classes of the vertical and horizontal subbundles of $T Z_{\pm}(M)$. We saw in the proof of Proposition~\ref{cohtwist4} that $Z_\pm(M)$ is the projectivization
of a~rank two complex vector bundle $E$ with total Chern class
\[
1+ \alpha - \frac{p_1 \pm 2\mathsf{e} - \alpha^2}{4}.
\]

As $T^{\rm v}(\mathbb{P}(E))\oplus \C \cong \mathcal{O}_{E^*}(1) \otimes p^*E$ (where $\mathbb{P}(E) \xrightarrow{p} M$ is the projection and $\mathcal{O}_{E^*}(1)$
is the canonical line bundle restricting to $\mathcal{O}(1)$ on each fiber~\cite[Definition~15.13]{D12}) we
have that
\[
c_1\big(T^{\rm v}(Z_\pm(M))\big) = c_1(E) + 2 c_1(\mathcal{O}_{E^*}(1)) = \alpha + 2x.
\]

As for the horizontal bundle, consider the pullback diagrams
\[
\begin{tikzcd}
Z_+(M) \ar[r] \ar[d] & Z_+\big({\mathcal B}{\rm Spin}^{\rm c}(4)\big) \ar[r,"\tau"] \ar[d] & Z_+({\mathcal B}{\rm SO}(4))={\mathcal B}{\rm U}(2)
\ar[d,"{\mathcal B}\iota"] \\
M \ar[r] & {\mathcal B}{\rm Spin}^{\rm c}(4) \ar[r,"{\mathcal B}\pi"] & {\mathcal B}{\rm SO}(4).
\end{tikzcd}
\]
Note that $Z_+({\mathcal B}{\rm SO}(4))= {\mathcal E}{\rm SO}(4) \times_{{\rm SO}(4)} {\rm SO}(4)/{\rm U}(2)={\mathcal B}{\rm U}(2)$, and under this identification the projection $Z_+({\mathcal B}{\rm SO}(4)) \to {\mathcal B}{\rm SO}(4)$ is the map ${\mathcal B}\iota$ induced by the inclusion ${\rm U}(2) \subset {\rm SO}(4)$. Moreover, the universal bundle ${\mathcal E}{\rm SO}(4) \times_{{\rm U}(2)} \R^4 \to {\mathcal E}{\rm SO}(4)/{\rm U}(2)$
pulls back to $T^{\rm h} Z_+(M)$ under the composite map $Z_+(M) \to Z_+\big({\mathcal B}{\rm Spin}^{\rm c}(4)\big) \xrightarrow{\tau} {\mathcal B}{\rm U}(2)$ (cf.\ Lemma~\ref{positchern}).

We saw in the proof of Proposition~\ref{cohtwist4} that
\[
H^*\big(Z_+\big({\mathcal B}{\rm Spin}^{\rm c}(4)\big)\big)\cong \Z[\alpha,S_1,S_2,x]/\big(x^2+\alpha x - S_2\big).
\]
In terms of this identification we have $\tau^*(c_2)={\mathcal B}\pi^*\mathsf{e} = -S_1+S_2$
and
\[
\tau^*(c_1) = \gamma x+\lambda \alpha \qquad \text{for some } \gamma, \lambda \in \Z.
\]
As the bundle classified by $\CP^1 \cong {\rm SO}(4)/{\rm U}(2) \to {\mathcal B}{\rm U}(2)$ has $c_1$ equal to two times the orientation class (see Lemma~\ref{positchern}), we have $\gamma=2$.
Since $p_1\in H^4({\mathcal B}{\rm SO}(4))$ maps to $c_1^2-2c_2$ in $H^4({\mathcal B}{\rm U}(2))$, and to $2S_1+2S_2 + \alpha^2$ in $H^4\big({\mathcal B}{\rm Spin}^{\rm c}(4)\big)$ we see that
\[
(2x+\lambda \alpha)^2 +2S_1-2S_2 = 2S_1+2S_2+ \alpha^2 \quad \Leftrightarrow \quad 4x^2 + 4\lambda \alpha x + \lambda^2 \alpha^2 = 4S_2+\alpha^2.
\]
Hence $\lambda=1$, yielding the formula \eqref{chernform} in the case of $Z_+(M)$.

Let $k \in {\rm O}(4)$ be an element with determinant $-1$ and let
$\varphi \colon {\rm SO}(4) \to {\rm SO}(4)$ denote conjugation by $k$. As a homogeneous space
of ${\rm SO}(4)$ we have $J_-\big(\R^4\big)={\rm SO}(4)/\varphi({\rm U}(2))$. Hence
$Z_-({\mathcal B}{\rm SO}(4))= {\mathcal E}{\rm SO}(4) \times_{{\rm SO}(4)}\varphi({\rm U}(2)) = {\mathcal B}(\varphi({\rm U}(2))$
and we have a commutative diagram
\[
\begin{tikzcd}
Z_-(M) \ar[r] \ar[d] & Z_-\big({\mathcal B}{\rm Spin}^{\rm c}(4)\big) \ar[r,"\tau'"] \ar[d] \ar[rr,bend left, "\tau"] & {\mathcal B}(\varphi({\rm U}(2)) \ar[r, "{\mathcal B}\varphi^{-1}"] \ar[d] & {\mathcal B}{\rm U}(2) \ar[d] \\
M \ar[r] & {\mathcal B}{\rm Spin}^{\rm c}(4) \ar[r, "{\mathcal B}\pi"] & {\mathcal B}{\rm SO}(4) \ar[r,"{\mathcal B}\varphi^{-1}"] & {\mathcal B}{\rm SO}(4)
\end{tikzcd}
\]
with the middle and left-hand squares both pullback squares
and the right hand square an
isomorphism of fiber bundles induced by the automorphism $\varphi^{-1}$.
Note that ${\mathcal B}\varphi^{-1}$ is covered by an isomorphism between the (complex) universal
bundles ${\mathcal E}{\rm SO}(4) \times_{\varphi({\rm U}(2))} \R^4 \to {\mathcal E}{\rm SO}(4)\times_{{\rm U}(2)} \R^4$,
so that the composition $Z_-(M) \to Z_-\big({\mathcal B}{\rm Spin}^{\rm c}(4)\big) \xrightarrow{\tau} {\mathcal B}{\rm U}(2)$
classifies $T^{\rm h} Z_-(M)$.

Considering the action of ${\mathcal B}\varphi$ on the maximal torus we see that ${\mathcal B}\varphi \colon
{\mathcal B}{\rm SO}(4) \to {\mathcal B}{\rm SO}(4)$ has the following effect on cohomology:
\[{\mathcal B}\varphi^*(p_i)=p_i, \qquad {\mathcal B}\varphi^* (\mathsf{e}) = -\mathsf{e}.\]

Arguing as before, we conclude that $\tau^*(c_2)=-{\mathcal B}\pi^*(\mathsf{e})=S_1-S_2$ and
$\tau^*(c_1)=2x+\alpha$,
leading to the formula~\eqref{chernform} in the case of~$Z_-(M)$.
\end{proof}

\begin{Remark}\label{hitchin} \quad
\begin{itemize}\itemsep=0pt

\item The formulas in Proposition~\ref{chernclasses} for $c_2$ and $c_3$ of $Z_-(M)$ differ by a sign from those given in~\cite[p.~135]{Hi81}; however, the fundamental class being used in loc.\ cit.\ seems to also differ by a sign from the one induced by the Atiyah--Hitchin--Singer almost complex structure on $Z_-(M)$. Hence the values of the Chern numbers given in~\cite[equation~(1.5)]{Hi81} coincide with those obtained with the Chern classes in Proposition~\ref{chernclasses}.

\item In~\cite[equation~(1.4)]{Hi81}, there is a description of the real cohomology ring of the negative twistor space of a four-manifold $X$, as the free $H^*(X;\R)$-module generated by $h=\frac{1}{2}c_1\big(T^{\rm v} Z_-(X)\big)$ subject only to the relation $h^2 = \tfrac{1}{2}\mathsf{e}(X) - \tfrac{1}{4}p_1(X)$.
However, this description is at odds with the identification of the negative twistor space of $\CP^2$ as the full flag variety ${\rm U}(3)/\left( {\rm U}(1)\times {\rm U}(1) \times {\rm U}(1) \right)$~\cite[p.~133]{Hi81},~\cite[p.~217, Example~1]{Sa84}. Namely, $H^*\big(\CP^2;\R\big) \cong \R[x]/\big(x^3\big)$, where $\big\langle x^2, \big[\CP^2\big] \big\rangle = 1$, so $\mathsf{e}\big(\CP^2\big) = p_1\big(\CP^2\big) = 3x^2$. Therefore the real cohomology of the negative twistor space is,
according to~\cite[equation~(1.4)]{Hi81}, $\R[x,h]/\big(x^3, h^2 - \tfrac{3}{4}x^2\big)$. One can check directly, however, that this (graded) ring is not isomorphic to the cohomology of the flag variety, i.e., to $\R[x,h]/\big(x^2+xh+h^2, x^2 h + h x^2\big)$.

In terms of Propositions~\ref{cohtwist4} and~\ref{chernclasses}, we have $h=x+\frac{\alpha}{2}$ and therefore $h^2=\frac{1}{4} p_1(x) - \frac{1}{2} \mathsf{e}(X)$.
\end{itemize}
\end{Remark}

\begin{Remark}
Let $M$ be an almost complex Riemannian six-manifold.
The fiberwise diffeomorphism $a\colon Z_+(M) \to Z_-(M)$ sending
$J_x \mapsto -J_x \in Z(T_xM)$ induces an anti-holomorphic map
between the fibers and hence
\[
a^*(c_1(T_v(Z_-(M)))) = -c_1(T_v(Z_+(M))).
\]
The spaces $Z_\pm(M)$ are projectivizations of rank 4 complex vector bundles with first Chern class $2c_1(M)$ (see Corollary~\ref{cohomtwis}), and so as in Proposition~\ref{chernclasses} we see that $c_1(T_v Z_\pm(M))=2c_1(M)+4x$. Therefore
\[
a^*(2c_1 + 4x) = -2c_1 - 4x \ \Rightarrow \ a^*(4x) = -4x-4c_1.
\]
As there is no torsion in $H^2$ of the universal example $Z({\mathcal B}{\rm U}(3))$, it follows that $a^*(x)=-x-c_1$. One can check that $a^*$ yields an isomorphism between the $H^*(M)$-algebras of Corollary~\ref{cohomtwis}.
\end{Remark}

In view of the previous remark, from now on we will restrict our attention to the positive twistor space of almost complex Riemannian six-manifolds.

\section[On the homotopy type of the space of almost complex structures on six-manifolds]{On the homotopy type of the space\\ of almost complex structures on six-manifolds}

In this section we study the homotopy type of the components of the space of almost complex structures on connected six-manifolds $M$, complementing our previous treatment~\cite{FGM21} of the case when $c_1=0$.

\begin{Theorem}\label{fibration}
Let $M$ be an oriented Riemannian manifold of dimension four or six, $J$ an orthogonal almost complex structure on $M$ compatible with the orientation, and $J(M)$ its component in the space of orthogonal almost complex structures. Let $S_{\C}^+(M)$ be the positive spinor bundle on $M$ determined by $J$. Then the map
\[
\Gamma\big( S\big(S_\C^+(M)\big)\big) \to J(M)
\]
is surjective, and for any lift $s$ of $J$,
\[
\Map\big(M,S^1\big) \to \Gamma\big( S\big(S_\C^+\big)\big)_s \to J(M)
\]
is a fiber sequence.
\end{Theorem}
\begin{proof}
Identifying $S_\C^+(M)$ with $\oplus_{k \text{ even}} \Lambda^k_\C(TM)$, Remark~\ref{lift} shows that
the constant section $s \colon M \to \Lambda^0_\C(TM) \subset S_\C^+(M)$ defined by
$s(x)=1$ lifts the section $J$. As $S\big(S_\C^+\big) \to Z_+(M)$ is a~fiberwise fibration over $M$, the map induced on sections is a~fibration. Since the fiber over $J$ is nonempty, any choice of lift of $J$ identifies the
fiber with $\Map\big(M,S^1\big)$.
\end{proof}

\begin{Remark} The above argument remains valid in all dimensions if one replaces $S\big(S_{\C}^+(M)\big)$ with $S\big(S_{\C}^+(M)\big) \cap {\mathcal P}S ^+_{\C}(M)$. \end{Remark}

Now let $(M,J)$ be a closed almost complex Riemannian six-manifold with $H^1(M;\ZZ) = 0$ (equivalently, $b_1 = 0$) satisfying $c_1c_2 - c_3 \neq 0$. By Theorem~\ref{fibration}, we have the fibration $\Map\big(M,S^1\big) \to \Gamma\big( S\big(S_\C^+(M)\big)\big)_s \to J(M)$. Since $H^1(M;\ZZ) = 0$, evaluation at a chosen basepoint gives a homotopy equivalence $\Map\big(M,S^1\big) \xrightarrow{{\rm ev}} S^1$, so we have a principal fiber sequence $S^1 \to \Gamma\big(S\big(S_\C^+(M)\big)\big)_s \to J(M)$.

We can thus consider instead the fiber sequence \begin{equation}\label{fibration2} \Gamma(S\big(S_{\C}^+\big))_s \to J(M) \to {\mathcal B}S ^1 \simeq \CP^\infty. \end{equation} Now, $S\big(S_{\C}^+\big)$ is an oriented fiber bundle over $M$ with fiber $S^7$, and hence it is classified by a map $M \to {\mathcal B}{\rm Aut}^+\big(S^7\big)$ to the classifying space of orientation-preserving homotopy automorphisms of~$S^7$. It is known that ${\mathcal B}{\rm Aut}^+\big(S^7\big)_{\QQ} \simeq K(\QQ, 8)$, where the subscript of $\QQ$ denotes (Sullivan) rationalization~\cite[Section~11]{S77}. By~\cite[Theorem~5.3]{M87}, the rationalization of $\Gamma\big(S\big(S_{\C}^+\big)\big)_s$ is homotopy equivalent to (the connected component corresponding to $s$ in) the space of sections of the fiberwise rationalized $S^7$ bundle over $M$. Since $H^8(M;\QQ) = 0$, this latter bundle is trivial, and hence $\big(\Gamma\big(S\big(S_{\C}^+\big)\big)_s\big)_{\QQ}$ is homotopy equivalent to a connected component of $\Map\big(M, S^7_{\QQ}\big)$. Denoting the Betti numbers of $M$ by $b_i$, by Thom's theorem on the space of maps into an Eilenberg--Maclane space, $\Map\big(M, S^7_{\QQ}\big)$ is in fact connected and has the homotopy type of $S^1_{\QQ} \times S^7_{\QQ} \times K(\QQ, 3)^{b_4} \times K(\QQ, 4)^{b_3} \times K(\QQ,5)^{b_2}$.

We can describe the fundamental group of $J(M)$ by using a theorem of Crabb and Sutherland:

\begin{Theorem}[{\cite[Theorem~2.12(i)]{CS84}}] Let $M$ be a closed connected $2n$-manifold and $\xi$ a complex $(n+1)$-plane bundle over $M$. Denote by $N\xi$ any component of the space of sections of $\mathbb{P}(\xi)$ whose elements lift to sections of $\xi$. Then $\pi_1(N\xi)$ is a central extension \[0 \to \ZZ/\big(\big\langle c_n(\xi), [M] \big\rangle \big) \to \pi_1(N\xi) \to H^1(M) \to 0.\]
\end{Theorem}

By Theorem~\ref{fibration}, we can apply this to $\xi = S_{\C}^+(M) \cong \C \oplus \Lambda_{\C}^2 TM$, which satisfies ${N \xi = J(M)}$. Therefore, the fundamental group of $J(M)$ is given by $\ZZ$ modulo $\int_M\! c_1c_2-c_3$ (where ${c_i \! =\! c_i(TM)}$), see Corollary~\ref{cohomtwis}. By~\cite{M87}, $J(M)$ is a nilpotent space as it is the space of sections of a fibration with nilpotent fiber \big(namely $\CP^3$\big) over a finite-dimensional base. We may thus rationalize the fibration (\ref{fibration2}) to obtain the fibration \[\big(\Gamma\big(S\big(S_{\C}^+\big)\big)_s\big)_{\Q} \to J(M)_{\Q} \to K(\Q, 2).\] Consider the degree one rational class corresponding to the factor $S^1_{\QQ}$ in the fiber. If $c_1c_2 - c_3 \neq 0$, then since $\pi_1(J(M)_\Q) = 0$, this class must hit (a nonzero multiple of) the degree two generator in $K(\Q,2)$ in the Serre spectral sequence. A model for the total space is given by the tensor product of the base and fiber's models with a perturbed differential on the fiber generators~\cite[Section~4]{S77}; the previous sentence thus tells us that a model for $J(M)$ is of the form
\[\big(\Lambda\big(x_2, z_1, z_7, z_3^{i}, z_4^{j}, z_5^{k}\big),\, {\rm d}z_1 = x_2, \, {\rm d}z_7 \in (x_2), \, {\rm d}z_3^i \in (x_2), \, {\rm d}z_4^j \in (x_2),\, {\rm d}z_5^k \in (x_2)\big),
\] where $(x_2)$ denotes the ideal generated by $x_2$, and $i$, $j$, $k$ range over sets of size $b_4$, $b_3$, $b_2$, respectively.

Notice that such a model is not minimal, and in fact a minimal model is obtained by quotienting out the differential ideal generated by~$z_1$. Indeed, by the argument in~\cite[Proposition~2]{VS76} we see that the map $(\Lambda, {\rm d}) \to \big(\Lambda/(z_1, {\rm d}z_1)\Lambda, \bar{{\rm d}}\big)$ induces an isomorphism on cohomology.
Here by $(\Lambda,{\rm d})$ we denote the differential graded algebra displayed above, and by $\bar{{\rm d}}$ the induced differential on the quotient. We see that
\[
(\Lambda/(z_1, x_2)\Lambda, \bar{{\rm d}}) \cong \big(\Lambda\big(z_7, z_3^{i}, z_4^{j}, z_5^{k}\big), \bar{{\rm d}}=0\big),
\] where the right-hand side is minimal. To summarize, we have the following:

\begin{Theorem}\label{Qhomotopy6} Let $M$ be a connected closed six-manifold with $b_1 = 0$ equipped with an almost complex structure $J$ such that $\int_M c_1c_2-c_3 \neq 0$. Then the space of almost complex structures on~$M$ in the component of $J$ is a nilpotent space with finite cyclic fundamental group and minimal model given by $\big(\Lambda\big(z_7, z_3^{i}, z_4^{j}, z_5^{k}\big), {\rm d} = 0\big)$; here $i$, $j$, $k$ range over sets of size $b_4(M)$, $b_3(M)$, $b_2(M)$ $(=b_4(M))$, respectively. In particular, this space is formal.
\end{Theorem}

We emphasize that the above makes no nilpotency assumption on $M$. Alternatively, a model for (a given component of) the space of almost complex structures can be obtained with the Haefliger--Sullivan model for the space of sections of a fibration~\cite[Section~11]{S77},~\cite{H82}, applied to $\CP^3 \to Z_+(M) \to M$. Our approach above lets one quickly describe the minimal model of~$J(M)$ using the geometry of the setup.

\begin{Example}\label{examplesQhomotopy} Using Theorem~\ref{Qhomotopy6}, we immediately obtain the following:
\begin{enumerate}\itemsep=0pt
\item Observe that the space of almost complex structures inducing a given orientation on the connected sum $g\big(S^3 \times S^3\big)$ is connected. By Theorem~\ref{Qhomotopy6}, for $g \neq 1$, the rationalization of this space is $K(\QQ, 7) \times K(\QQ, 4)^{2g}$. In particular, the space of almost complex structures on $S^6$ has the same rational homotopy type as $S^7$, i.e., $K(\QQ,7)$; cf.~\cite{FGM21} where it is shown that a certain natural inclusion $\RP^7 \hookrightarrow J\big(S^6\big)$ induces an isomorphism on rational homotopy groups and fundamental groups. (Note that the covering $S^7 \to \RP^7$ induces an isomorphism on rational homotopy groups.)

For $g=1$, one can calculate directly using the Haefliger--Sullivan model that the space of almost complex structures on $S^3 \times S^3$ has the same rational homotopy type as $S^1 \times \CP^3 \times (\HP^\infty)^2$, cf.~\cite[Example~5.4]{FGM21}.

\item The components of almost complex structures on $\mathbb{CP}^3$ are parametrized by $c_1 = 2kx$, $c_2 = \big(2k^2-2\big)x^2$, where $x$ is the generator of $H^2\big(\CP^3\big)$ such that $\big\langle x^3, \big[\CP^3\big]\big\rangle = 1$ (i.e., $x = c_1(\mathcal{O}(1))$ for the standard complex structure). Hence $c_1c_2 - c_3 = 4k^3 - 4k - 4$, which is never zero. Hence every component of almost complex structures has rationalization $K(\QQ,7) \times K(\QQ,5) \times K(\QQ,3)$. Contrast this with the trivial $\CP^3$ bundle over $\CP^3$, whose (infinitely many) components of sections exhibit two distinct rational homotopy types~\cite[Example~3.4]{MR85}; the space of maps $\CP^3 \to \CP^3$ homotopic to any fixed essential map has rationalization $K(\QQ,7)\times K(\QQ,5)\times K(\QQ,3)$. \end{enumerate} \end{Example}

\section[Intersections of almost complex structures on six-manifolds as sections of the twistor space]{Intersections of almost complex structures on six-manifolds\\ as sections of the twistor space}

We will now use the results of Section~\ref{section3} to calculate the homological intersection of two almost complex structures on a given six-manifold $M$, after identifying them with the image of $M$ under the corresponding section of the twistor bundle.

\begin{Lemma}\label{poincaredual}
Let $E \xrightarrow{\pi} M$ be a complex vector bundle of rank $n$ over a closed oriented mani\-fold~$M$. Let $s \colon M \to \mathbb{P}(E \oplus \C)$ denote the canonical section given pointwise by $[0: \cdots : 0 : 1]$.
Then the Poincar\'e dual of $s_*[M]$ in $H^*(\mathbb{P}(E\oplus \C)) \cong H^*(M)[x]/\big(x^{n+1} + c_1 x^n + \cdots + c_n x\big)$ is
\[
x^n + c_1x^{n-1} + \cdots + c_n,
\]
where $c_i$ denote the Chern classes of $E$.\footnote{Here the orientation on $\mathbb{P}(E\oplus \C)$ is understood to be the one induced by the orientation on $M$ and the canonical orientation on the fiber $\CP^n$ corresponding to $\langle x^n, \big[\CP^n\big] \rangle = 1$.}
\end{Lemma}
\begin{proof}
The Poincar\'e dual to $s_*[M]$ is the image of the Thom class of the normal bundle to the submanifold $s(M) \subset \mathbb{P}(E \oplus \C)$ under the first map below in the long exact sequence of a pair,
\[
H^{2n}(\mathbb{P}(E \oplus \C), \mathbb{P}(E \oplus \C) \setminus s(M)) \to H^{2n}(\mathbb{P}(E \oplus \C)) \to H^{2n} (\mathbb{P}(E \oplus \C) \setminus s(M)),
\]
where $H^{2n}(\mathbb{P}(E \oplus \C), \mathbb{P}(E \oplus \C) \setminus s(M))$ has been identified with $H^{2n}(U, U \setminus s(M))$ for some tubular neighborhood $U$ of $s(M)$ via excision; see~\cite[Section~6.11]{Br}. This is a nonzero class, as~$s_*[M]$ is nonzero since it projects to~$[M]$.

Now notice that $\mathbb{P}(E \oplus \C) \setminus s(M)$ deformation retracts to $\mathbb{P}(E)$. Therefore the effect of the second map on cohomology is the natural quotient
\[ H^{2n}(M)[x]/\big(x^{n+1} + c_1 x^n + \cdots + c_n x\big) \to H^{2n}(M)[x]/\big(x^n + c_1 x^{n-1} + \cdots + c_n\big), \]
and thus the kernel is generated by $x^n + c_1x^{n-1} + \cdots + c_n$. Hence the Thom class maps to $x^n + c_1x^{n-1} + \cdots + c_n$ up to sign. To check the sign, we restrict to a fiber and evaluate against the fundamental class $\big[\CP^n\big]$. Since the Thom class is the orientation of the normal bundle, which is given by the orientation of the fibers, and $\big\langle x^n, \big[\CP^n\big] \big\rangle = 1$, we see that the sign is positive. \end{proof}

\begin{Theorem}\label{selfintersection}
Let $s$ be a section of the positive twistor bundle over a closed Riemannian six-manifold $M$, and give $M$ the orientation determined by $s$. Then the homological self-intersection number $s_*[M] \cdot s_*[M]$ of $s(M)$ in $Z_+(M)$ is given by $\int_M c_1c_2 - c_3$, where the $c_i$ are the Chern classes of the almost complex structure determined by $s$.

\end{Theorem}
\begin{proof} In terms of the identification in Proposition~\ref{identification}, $s$ is the canonical section of $\mathbb{P}\big(S_{\C}^+\big) = \mathbb{P}\big(\Lambda_\C^2 TM \oplus \C\big)$ (see the proof of Theorem~\ref{fibration}), where $TM$ is equipped with the complex structure determined by $s$. To compute the homological self-intersection of $s(M)$, we integrate the square of the expression obtained in Lemma~\ref{poincaredual} for its Poincar\'e dual in $Z_+(M)$. We use Corollary~\ref{cohomtwis}, and observe that since $x\big(x^3 + c_1\big(S_{\C}^+\big)x^2 + c_2\big(S_{\C}^+\big)x + c_3\big(S_{\C}^+\big)\big) = 0$ in $H^*(Z_+(M))$, we have
\begin{align*}
\big(x^3 + c_1\big(S_{\C}^+\big)x^2 + c_2\big(S_{\C}^+\big)x + c_3\big(S_{\C}^+\big)\big)^2 &= c_3\big(S_{\C}^+\big)\big(x^3 + c_1\big(S_{\C}^+\big)x^2 + c_2\big(S_{\C}^+\big)x + c_3\big(S_{\C}^+\big)\big) \\ &= c_3\big(S_{\C}^+\big) x^3 = (c_1c_2 - c_3)x^3. \tag*{\qed}
\end{align*}\renewcommand{\qed}{}
\end{proof}

Using Theorem~\ref{selfintersection}, we obtain the following:

\begin{Example}\label{intersectiontwistor4d} \quad\samepage
\begin{enumerate}\itemsep=0pt \item
For any closed oriented Riemannian four-manifold $M$, by Proposition~\ref{chernclasses}, we have the following Chern numbers of the Atiyah--Hitchin--Singer almost complex structure on $Z_+(M)$:
\begin{gather*}
\int_{Z_+(M)} c_1(Z_+(M))c_2(Z_+(M)) = \int_{Z_+(M)} (4x+2\alpha)(p_1+3\mathsf{e}) \\
\hphantom{\int_{Z_+(M)} c_1(Z_+(M))c_2(Z_+(M))}{}
= \left( \int_{\CP^1}4x \right)\left( \int_{M} (p_1(M)+3\mathsf{e}(M)) \right) \\
\hphantom{\int_{Z_+(M)} c_1(Z_+(M))c_2(Z_+(M))}{}
 = 12(\sigma(M)+\chi(M)), \\
\int_{Z_+(M)} c_3(Z_+(M)) = \int_{Z_+(M)} (2x+\alpha)\mathsf{e} = 2 \chi(M).
\end{gather*}

Hence the self-intersection of the Atiyah--Hitchin--Singer almost complex structure on $Z_+(M)$ in $Z_+(Z_+(M))$ is $12\sigma(M) + 10\chi(M)$.

By~\cite[Lemma~3]{Ar97}, closed oriented four-manifolds admitting a maximally nonintegrable almost complex structure satisfy $5\chi + 6\sigma = 0$. Hence for these manifolds, the corresponding self-intersection number is zero.

\item By a theorem of Taubes~\cite[Theorem~1.1]{T92}, given any closed oriented four-manifold $M$, for every large enough $k$ the manifold $N = \overline{M} \# k \CP^2$ carries a self-dual metric. Hence its negative twistor space $X = Z_-(N)$ is a complex manifold~\cite[Theorem~4.1]{AHS78}. It satisfies \begin{align*} \int_X c_1(X)c_2(X) = 12(\chi(N) - \sigma(N)), \qquad \int_X c_3(X) = 2\chi(N).\end{align*}

Indeed, we obtain this again from Proposition~\ref{chernclasses}, namely \begin{align*} \int_X c_1(X)c_2(X) &= \int_X 2(\alpha +2x)(p_1(N)-3\mathsf{e}(N)) \\ &= \left( \int_{\CP^1}4x \right) \left( \int_{\overline{N}} (p_1(N)-3\mathsf{e}(N)) \right) = -12(\sigma(N)-\chi(N)).\end{align*} This agrees with~\cite[equation~(1.5)]{Hi81}. Alternatively, from~\cite[Section~4]{EaSi93} we see that the holomorphic Euler characteristic of $X$ is given by $1 - b_1(N) + b_2^+(N)$. Since $\chi(N) = 2+b_2^+(N) +b_2^-(N) - 2b_1(N)$ and $\sigma(N) = b_2^+(N) - b_2^-(N)$, we have $\tfrac{1}{2}(\chi(N) - \sigma(N)) = 1 - b_1+b_2^-$, which equals the Todd genus $\int_X \tfrac{1}{24}c_1c_2(X)$. The Euler characteristic of $X$ is twice that of~$N$, as seen immediately, e.g., using the Leray--Hirsch theorem.

From here we obtain \[\int_X c_1(X)c_2(X) - c_3(X) = 10\chi(M) + 12\sigma(M) - 2k.\] Hence for any $M$, for large enough $k$ we obtain compact complex threefolds $X$ with negative homological self-intersection in their own positive twistor space.
\end{enumerate}
\end{Example}

\begin{Example}\label{exampleselfintersection}\quad
\begin{enumerate}\itemsep=0pt
\item For the connected sum $g\big(S^3 \times S^3\big)$, the self-intersection number of any almost complex structure in $Z_+\big(g\big(S^3 \times S^3\big)\big)$ is $2g-2$. In particular, for $S^6$ we have $-2$.
\item $S^2 \times S^4$ admits a unique homotopy class of almost complex structures for every choice of $c_1 = 2k\alpha$, where $\alpha$ generates $H^2$. Since $p_1 = 0$ we see that $c_2 = 0$ for any $k \in \Z$, and so the self-intersection is $-4$ for any homotopy class of almost complex structure. Likewise, for $g\big( S^2 \times S^4 \big)$, the self-intersection of any almost complex structure is $-2g-2$.
\item Recall that the homotopy classes of almost complex structures on $\mathbb{CP}^3$ are parametrized by $c_1 = 2kx$, $c_2 = \big(2k^2-2\big)x^2$. The self-intersection numbers are given by $4k^3 - 4k - 4$; in particular, for the standard complex structure it is 20.

\item To obtain more examples of compact complex threefolds with negative self-intersection, one can blow up any given compact complex threefold $X$ at a point sufficiently many times. Indeed, the blowup at one point, diffeomorphic to $X \# \overline{\CP^3}$, satisfies
\[ \int_{X \# \overline{\CP^3}} c_1c_2\big(X \# \overline{\CP^3}\big) = \int_X c_1c_2(X)\] by invariance of the Todd genus under bimeromorphisms, and the Euler characteristic is two larger than that of $X$. Note that if $X$ is K\"ahler, this gives examples of K\"ahler threefolds with negative self-intersection.
\item In~\cite{Le99}, LeBrun builds examples of compact complex threefolds $X_m$, $m>0$, each diffeomorphic to $K3\times S^2$, with $\int_{X_m} c_1c_2(X_m) = 48m$. Hence the homological self-intersection number is given by $48(m-1)$.
\item By~\cite[Corollary 3.1]{GS14}, for a connected symplectic six-manifold $M$ with a Hamiltonian circle action with isolated fixed points, we have $\int_M c_1c_2 = 24$ and so $\int_M c_1c_2 - c_3 = 24-\chi(M)$. In particular, this formula applies to toric six-manifolds.

\end{enumerate} \end{Example}

Recall the following theorem of Michelsohn, Salamon, Atiyah--Hitchin--Singer:

\begin{Theorem}[{\cite[Theorem~9.11]{LM}\label{sectiontheorem},~\cite[Theorem~1]{Sa95}}] A section $s$ of the positive twistor space $Z_+(M) \to M$ is pseudoholomorphic, with respect to the almost complex structure $J$ corresponding to $s$ on $M$ and the Atiyah--Hitchin--Singer almost complex structure on $Z_+(M)$, if and only if $J$ is integrable. \end{Theorem}

Now suppose we have an integrable orthogonal $J$ on a Riemannian six-manifold $M$. Then $s$ gives us an embedding of $M$ into $Z_+(M)$ as an almost complex submanifold. If we were able to perturb $s(M)$ to another almost complex submanifold that intersects $s(M)$ transversally, then the homological intersection number $s_*[M] \cdot s_*[M]$ would have to be $\geq 0$. A way to obtain perturbations of $s(M)$ to almost complex submanifolds is by translating $J$ by an isometry of $(M,g)$, but any perturbation to a transverse representative would do. Note that it is unique to dimension six that $s_*[M] \cdot s_*[M]$ inside $Z_+(M)$ gives an integer.

\begin{Corollary}\label{deform} Suppose a closed Riemannian six-manifold $M$ carries an integrable complex structure, corresponding to a section $s$ of $Z_+(M)$, with $\int_M c_1c_2 - c_3 < 0$. Then there is no representative of the homology class $s_*[M]$ given by an almost complex submanifold of $Z_+(M)$ which intersects $s(M)$ transversally. \end{Corollary}

\begin{Example}
Consider $S^6$ with the round metric. Take the canonical almost complex structure induced by the octonions, and consider its orbit under the action of the isometry group ${\rm SO}(7)$, diffeomorphic to $\RP^7$. After identifying the twistor space of round $S^6$ with the Grassmannian of oriented 2-planes in $\R^8$, Calabi and Gluck~\cite{CG} show that these ``octonion'' $J$'s send $S^6$ to the family of 2-planes of the form $(o, v)$, where $o$ is a fixed unit octonion depending on $J$ (so, the image of $S^6$ is the unit sphere in the hyperplane orthogonal to a given line in $\R^8$.) Two such nearby families intersect transversally. Therefore, the $-2$ homological self-intersection obtained in Example~\ref{exampleselfintersection} shows the octonion $J$'s are not integrable. Of course, this also follows by a direct computation of the Nijenhuis tensor, and the result is a special case of the classical theorem of Blanchard~\cite{Bl53} and LeBrun~\cite{Le87} that there is no integrable complex structure on $S^6$ orthogonal with respect to the round metric.
\end{Example}

In order to study intersections of almost complex structures in different components of the space of sections of the twistor space we will now analyze the dependence of our formula for the cohomology (Corollary~\ref{cohomtwis}) on the choice of almost complex structure.

Consider the map of fibrations
\begin{equation}
\label{comparison}
\begin{tikzcd}
\CP^3 \ar[r,"\jmath"] \ar[d] & \CP^\infty \ar[d] \\
{\mathcal B}{\rm U}(3) \ar[r,"{\mathcal B}\kappa"] \ar[d] & {\mathcal B}{\rm Spin}^{\rm c}(6)\ar[d] \\
{\mathcal B}{\rm SO}(6) \ar[r,"="] & {\mathcal B}{\rm SO}(6).
\end{tikzcd}
\end{equation}
As the canonical map $\kappa \colon {\rm U}(3) \to \Spin^{\rm c}(6)$ is an isomorphism on $\pi_1$, we see
that $\jmath$ is an isomorphism on $\pi_2$ and hence the $5$-lemma implies that ${\mathcal B}\kappa$ is a
$7$-equivalence. Therefore the set of isomorphism classes of almost complex structures on six-manifolds is in natural bijective correspondence with the set of isomorphism classes of ${\rm spin}^{\rm c}$ structures via the canonical map from the former to the latter.

The set of isomorphism classes of ${\rm spin}^{\rm c}$-structures on a manifold $M$ is a torsor over the group $H^2(M)$
of isomorphism classes of complex line bundles
(corresponding to the action of the fiber on
the principal fibration along the right-hand column of \eqref{comparison}). Letting
$[s]$ denote the ${\rm spin}^{\rm c}$ structure associated to the almost complex structure $s\colon
M \to Z_+(M)$, and writing $a \cdot [s]$ for the action of the cohomology class
$a \in H^2(M)$ on $[s]$, we have the following result:

\begin{Lemma}
Let $M$ be an oriented Riemannian six-manifold and let $s,s'\colon M \to Z_+(M)$ be
two almost complex structures on $M$ compatible with the orientation of $M$.
Let $a \in H^2(M)$ be $($the unique cohomology class$)$ such that $[s']=a\cdot [s]$.
Then the canonical isomorphism between the expressions for $H^*(Z_+(M))$ given in Corollary~{\rm \ref{cohomtwis}} in terms of $s$, $s'$ is given by
\[
x' \mapsto x - a.
\]
\end{Lemma}
\begin{proof}
Let $E$ and $E'$ be the positive spinor bundles associated to the ${\rm spin}^{\rm c}$ structures~$[s]$ and~$[s']$, respectively, and let $L$ denote a complex line bundle with $c_1(L)=a$.

The relation $a \cdot [s]=[s']$ implies that there is an isomorphism
\[\psi \colon \ L \otimes E \to E'.\]
The kernel in $T_x M \otimes \C$ of right multiplication by a positive spinor
$v \in E_x$ is the same as the kernel of right multiplication by $\lambda \otimes v \in
L_x \otimes E_x$ for any $\lambda \neq 0$, so the canonical isomorphism given by the
composition
\[
\mathbb{P}(L\otimes E) \xrightarrow{\phi'} Z_+(M) \xrightarrow{\phi^{-1}} \mathbb{P}(E)
\]
is the obvious one sending the line $L \otimes (\C v)$ to $\C v$.

We conclude that on $Z_+(M)$ the tautological (Hopf) line bundles $H$ and $H'$ coming
respectively from the identifications with $\mathbb{P}(E)$ and $\mathbb{P}(E')$ satisfy
the relation
\[
H' \cong H \otimes L.
\]
Since $x'=c_1(H'^*)$ and $x=c_1(H^*)$, this completes the proof.
\end{proof}

Given a ${\rm spin}^{\rm c}$ structure $[s]$ we will write $c_1([s])$ for the Chern class of
the line bundle classified by the projection ${\mathcal B}{\rm Spin}^{\rm c}(2n) \xrightarrow{{\mathcal B}\delta} {\mathcal B}S ^1$.
Note that $c_1( a\cdot [s]) = 2a + c_1([s])$ so the difference class $a \in H^2(M)$
between two almost complex structures $s$, $s'$ is a specific ``square root'' of the
difference $c_1'-c_1$. We can suggestively write
\begin{equation}
\label{div2}
x' = x + \frac{c_1-c_1'}{2},
\end{equation}
which is unambiguous on rational cohomology (or when there is no $2$-torsion in $H^2(M)$).

\begin{Remark}
If we use $s\colon M \to Z_+(M)$ to write
\[
H^*(Z_+(M)) = H^*(M)[x]/\big(x^4 + 2c_1 x^3 + \big(c_1^2+c_2\big)x^2 + (c_1c_2-c_3) x\big),
\]
then $s^*(x)=0$. This follows from Remark~\ref{lift} which implies that
the pullback by $s$ of the Hopf bundle over $\mathbb P(E)$ is trivial. Using the isomorphism $H^2(Z_+(M)) \xrightarrow{\iota^*\oplus s^*} H^2\big(\CP^3\big) \oplus H^2(M)$ one sees that the first Chern class of the horizontal
plane bundle on $Z_+(M)$ is\footnote{This computation also follows from the arguments used in the proof of Proposition~\ref{chernclasses}.} $c_1(M)+2x$ (for any choice of almost complex structure) leading again to~\eqref{div2}.
\end{Remark}

\begin{Theorem}\label{intersection} Let $J$ and $J'$ be two orthogonal almost complex structures inducing the same orientation on a closed Riemannian six-manifold $M$, corresponding to sections $s$ and $s'$ of the twistor bundle. Denote by $c_i$, $c_i'$ their respective Chern classes. Then the intersection number $s_*[M] \cdot s'_*[M]$ in $Z_+(M)$ is given by \[\int_M \tfrac{1}{8}\big( c_1^3 + c_1^2 c_1' - c_1c_1'^2 - c_1'^3 \big) + \tfrac{1}{2} \big( c_1 c_2' + c_1' c_2' \big) - c_3.\]
\end{Theorem}

\begin{proof}
By Corollary~\ref{cohomtwis}, the presentation of $H^*(Z_+(M))$ with reference to $J$ is given by \[H^*(M)[x]/\big(x^4 + 2c_1x^3 + \big(c_1^2+c_2\big)x^2 + \big(c_1c_2-c_3\big)x\big),\] and with reference to $J'$ it is given by \[H^*(M)[x']/\big((x')^4 + 2c_1'(x')^3 + ((c_1')^2 + c_2')(x')^2 + (c_1'c_2'-c_3')x'\big).\] Combining Corollary~\ref{cohomtwis} with Lemma~\ref{poincaredual}, we have that the Poincar\'e dual of $s_*[M]$ is given by \[x^3 + 2c_1x^2 + \big(c_1^2+c_2\big)x + (c_1c_2-c_3),\] and the Poincar\'e dual of $s'_*[M]$ is given by \[(x')^3 + 2c_1'(x')^2 + \big((c_1')^2 + c_2'\big)x' + c_1'c_2'-c_3'.\] Recall from \eqref{div2} that $x' = x + \tfrac{1}{2}(c_1 - c_1')$. (Our calculation will be insensitive to torsion, so we may in fact work in the rational cohomology ring where there is no ambiguity in writing $\tfrac{1}{2}(c_1 - c_1')$.) Using this to rewrite the expression for $s'_*[M]$, and multiplying with $s_*[M]$, yields \[x^3\big( \tfrac{1}{8}\big( c_1^3 + c_1^2 c_1' - c_1c_1'^2 - c_1'^3 \big) + \tfrac{1}{2} ( c_1 c_2' + c_1' c_2' ) - c_3' \big).\] Pairing with the fundamental class of $Z_+(M)$ gives the result. \end{proof}

\begin{Remark} The expression for the intersection number in Theorem~\ref{intersection} is symmetric in $J$ and $J'$ as expected, which can be verified by using $c_1^2 - 2c_2 = p_1 = c_1'^2 - 2c_2'$ and $c_3 = c_3'$. \end{Remark}

\begin{Example} Let $J$ and $J'$ be two almost complex structures on $\CP^3$, with Chern classes $c_i$ and $c_i'$. Then $c_1 = 2kx$, $c_2 = \big(2k^2 - 2\big)x^2$, $c_3 = 4x^3$ and $c_1' = 2\ell x$, $c_2' = \big(2\ell^2 - 2\big)x^2$, $c_3' = 4x^3$, where $k$, $\ell$ are integers. From Theorem~\ref{intersection} we get that $J$ and $J'$, thought of as sections of the twistor bundle, have homological intersection given by \[k^3 + k^2 \ell + k \ell^2 + \ell^3 - 2k - 2\ell - 4.\] In particular, taking $J$ to be the standard complex structure, lying in the component determined by $k = 2$, we have that the homological intersection with any $J'$ with $\ell < 0$ is negative.
\end{Example}

\subsection*{Acknowledgements}
We thank Luis Fernandez and Scott Wilson for numerous illuminating discussions, and the referees for helpful comments. The first author was partially supported by FCT/Portugal through CAMGSD, IST-ID, projects UIDB/04459/2020 and UIDP/04459/2020. The second author would like to thank the Max Planck Institute for Mathematics in Bonn for its support, along with the Mittag-Leffler Institute in Djursholm for its hospitality during a visit to the ``Higher algebraic structures in algebra, topology and geometry'' program, where part of this work was carried out.

\pdfbookmark[1]{References}{ref}
\LastPageEnding

\end{document}